\documentclass[12pt,a4paper]{article}
\usepackage{theorem}%
\usepackage{graphicx}%
\newtheorem{lemma}{Lemma}

\newtheorem{theorem}[lemma]{Theorem}
\newtheorem{corollary}[lemma]{Corollary}
{\theorembodyfont{\upshape}}
{\theorembodyfont{\upshape}}
{\theorembodyfont{\upshape}}
{\theorembodyfont{\upshape}}
{\theorembodyfont{\upshape}}

\newcommand{\Z}{{\bf Z}}
\newcommand{\R}{{\bf R}}
\newcommand{\C}{{\bf C}}

\newcommand{\rme}{{\rm e}}

\newcommand{\sig}{\sigma}
\newcommand{\alp}{\alpha}

\newcommand{\gam}{\gamma}

\newcommand{\lam}{\lambda}
\newcommand{\del}{\delta}
\newcommand{\eps}{\varepsilon}

\newcommand{\Lam}{\Lambda}
\newcommand{\Dom}{{\rm Dom}}

\newcommand{\Spec}{{\rm Spec}}

\newcommand{\Ker}{{\rm Ker}}
\newcommand{\Ran}{{\rm Ran}}

\newcommand{\norm}{\Vert}

\newcommand{\Proof}{\underbar{Proof}{\hskip 0.1in}}

\usepackage{hyperref}
\title{An Indefinite Convection-Diffusion Operator}
\author{E B Davies}
\date{5 February 2007}
\begin{document}
\maketitle

\begin{abstract}
We give a mathematically rigorous analysis which confirms
the surprising results in a recent paper \cite{BOS} of
Benilov, O'Brien and Sazonov about the spectrum of a
highly singular non-self-adjoint operator that arises in a
problem in fluid mechanics.

MSC-class: 34Lxx, 76Rxx, 65F15, 65Q05

keywords: spectrum, non-self-adjoint, fluid mechanics,
pseudospectra, eigenvalue, basis.
\end{abstract}

\section{Introduction}

In a recent paper \cite{BOS} Benilov, O'Brien and Sazonov have shown
that the equation
\[
\frac{\partial f}{\partial t}= \eps\frac{\partial }{\partial
\theta}\left(\sin(\theta)\frac{\partial f}{\partial
\theta}\right)+\frac{\partial f}{\partial \theta}
\]
approximates the evolution of a liquid film inside a rotating
horizontal cylinder. The variable $\theta$ is taken to lie in
$[-\pi,\pi]$ and one assumes that the solutions $f$ are sufficiently
smooth and satisfy periodic boundary conditions.

The operator $H$ is highly non-self-adjoint (NSA) and it is not
amenable to standard elliptic techniques because the second order
coefficient is indefinite. For $\theta\in(0,\pi)$ the second order
term has a diffusive effect on the evolution but for $\theta\in
(-\pi,0)$ its effect is anti-diffusive. Many of the calculations in
\cite{BOS} are based on an asymptotic or WKB analysis for small
$\eps>0$, but this has dangers because infinite order approximate
eigenvalues of NSA operators need not be close to true eigenvalues.
Nor need eigenvalues computed by truncations of a highly
non-self-adjoint operator to large finite-dimensional subspaces by
standard methods be close to the eigenvalues of the original
operator; see \cite{BS,LOTS,TE} for examples and discussions of
their relationship to pseudospectra. Our goal in this paper is to
rederive some of the results in \cite{BOS} for a fixed positive
value of $\eps$ by a rigorous and non-asymptotic technique. We also
provide strong numerical evidence that the eigenvectors do not form
a basis. In our numerical calculations we take $\eps=0.1\,$.

Before proceeding we mention that closely related operators were
discussed in \cite[pp. 124-125, 406-408]{TE} and \cite{Ben}.

%%%%%%%%%%%%%%%%%%%%%%%%%%%%%%%%%%%%%%%%%%%%%%%%%%%%%%%%%%%%%%%%
\section{A Reformulation of the Problem}

We focus attention on the spectral properties of the operator
\[
(Hf)(\theta):=\eps\frac{\partial }{\partial
\theta}\left(\sin(\theta)\frac{\partial f}{\partial
\theta}\right)+\frac{\partial f}{\partial \theta}
\]
defined on all $C^2$ periodic functions $f\in L^2(-\pi,\pi)$. We
normally assume that $0<\eps<2$ for reasons explained in
Corollary~\ref{whysmalleps}. According to the WKB analysis of
\cite{BOS} the eigenvalue equation $ -iHf=\lam f $ has a sequence of
real eigenvalues which converge to the integers as $\eps\to 0$. This
suggests that the evolution equation is neutrally stable, but
Benilov et al. show that it exhibits explosive disturbances. This is
closely related to the pseudospectra of the operator. Our goal is to
prove that there are indeed real eigenvalues $\lam$ without
depending on WKB analysis, and to provide a simple and rigorous
method for computing them.

By expanding $f\in L^2(-\pi,\pi)$ in the form
\[
f(\theta)=\frac{1}{\sqrt{2\pi}} \sum_{n\in \Z} v_n\rme^{in\theta},
\]
one may rewrite the eigenvalue problem in the form $Av=\lam v$,
where $A=-iH$ is given by
\[
(Av)_n=\frac{\eps}{2} n(n-1)v_{n-1}-\frac{\eps}{2}
n(n+1)v_{n+1}+nv_n.
\]
The (unbounded) tridiagonal matrix $A$ is of the form
\[
A=\left(\begin{array}{ccc} A_-&0&0\\
0&0&0\\
0&0&A_+
\end{array}
\right)
\]
where $A_-$ acts in $l^2(\Z_-)$, the central $0$ acts in $\C$ and
$A_+$ acts in $l^2(\Z_+)$. The coefficient map $n\to -n$ induces a
unitary equivalence between $A_-$ and $-A_+$, so we only need study
the spectrum of $A_+$. Since $A_+^\ast=DA_+D^{-1}$ where
$D_{r,s}=\del_{r,s}(-1)^r$, $A_+$ and $A_+^\ast$ have the same
spectrum. We assume that $A_+$ has its natural maximal domain, and
observe that it is a closed operator. We will see that its
eigenvectors decrease more rapidly as $n\to +\infty$ the smaller
$\eps>0 $ is. We prove that the spectrum is discrete, i.e. that it
consists only of isolated eigenvalues of finite multiplicity, in
Section~\ref{SectResolvent}.

Benilov et al. correctly state in \cite{BOS} that one obtains very
poor numerical results if one simply truncates $A$ to produce a
finite matrix whose eigenvalues are then computed. We study the
matrix $A$ in a completely different manner.

The eigenvalue equation for $A_+$ may be written in the form
\begin{equation}
n(n-1)v_{n-1}-n(n+1)v_{n+1}+2\frac{n-\lam}{\eps}v_n=0.
\label{recurrence}
\end{equation}
The reality of the coefficients of (\ref{recurrence}) implies that
if $\lam$ is an eigenvalue then so is $\overline{\lam}$. It does
not, however, imply that the eigenvalues are real.

We confine attention to the solutions of (\ref{recurrence}) with
support in $\Z^+$, and regard the $n=1$ case, namely $\eps
v_2=(1-\lam) v_1$, as an initial condition. Since it is a second
order recurrence equation, the solution space of (\ref{recurrence})
is two-dimensional. We will see that one solution, often called the
subordinate solution, lies in $l^2(\Z_+)$, but no others do so if
$0<\eps<2$. We say that $\lam>0$ is an eigenvalue of $A_+$ if the
subordinate solution of the recurrence equation satisfies the
initial condition.

If one assumes that (\ref{recurrence}) has a solution of the form
$v_n=n^a(1+b/n+O(1/n^{2}))$, then one finds that $a=-1+1/\eps$ and
$b=\lam/\eps$. This motivates our next two lemmas. In the following
calculations we introduce constants $N_{\lam,\eps}^{(i)}$, and will
use the fact that they can always be increased without affecting the
results.

\begin{lemma}\label{GrowUpper} If $\lam\geq 0$, there exists $N =N_{\lam,\eps}^{(1)}$
such that if $v_n$ is a solution of (\ref{recurrence}) satisfying
$0<v_{N-i}\leq(N-i)^{a}$ for $i=1,\, 2$ where $a=-1+1/\eps$, then
$0<v_n\leq n^{a}$ for all $n\geq N$.
\end{lemma}

\Proof Suppose that $n\geq \lam+3$ and $0<v_{n-i}\leq (n-i)^{a}$ for
$i=1,\, 2$. Then
\begin{eqnarray*}
0\, <\, n^{-a}v_n &=& n^{-a}\left( \frac{n-2}{n}
v_{n-2}+2\frac{n-1-\lam}{\eps
n(n-1)}v_{n-1}\right)\\
&\leq & \left( 1-\frac{2}{n}\right)^{a+1}+%
2\frac{n-1-\lam}{\eps n(n-1)}\left( 1-\frac{1}{n}\right)^a\\
&= & 1-\frac{2\lam}{\eps n^2}+O(n^{-3})\\
&\leq& 1
\end{eqnarray*}
for all $n\geq N=N_{\lam,\eps}^{(1)}$. It follows inductively that
$0<v_n\leq n^a$ for all $n\geq N$.

\begin{corollary}\label{whysmalleps} If $\lam\geq 0$ and $\eps>2$ then every
solution of (\ref{recurrence}) lies in $l^2(\Z_+)$. In particular
every such $\lam$ is an eigenvalue of $A_+$.
\end{corollary}

\Proof Let $N=N_{\lam,\eps}^{(1)}$. Let $u$ be the solution of
(\ref{recurrence}) such that $u_{N-2}=0$ and $u_{N-1}=(N-1)^a$, and
let $v$ be the solution such that $v_{N-2}=(N-2)^a$ and $v_{N-1}=0$.
Since $a<-1/2$, Lemma~\ref{GrowUpper} implies that both lie in
$l^2(\Z_+)$. The space of all solutions is two-dimensional, so every
solution lies in $l^2(\Z_+)$, and this applies in particular to the
solution that satisfies the initial condition.

It is highly probably that one could avoid the above conclusion by
imposing boundary conditions at $+\infty$ if $\eps >2$, i.e. by
reducing the domain of $A_+$. We do not pursue this possibility.

\begin{lemma}\label{GrowLower} For every $\lam\in\R$ there
exists $N =N_{\lam,\eps}^{(2)}$ such that if $v_n$ is a solution of
(\ref{recurrence}) satisfying
\begin{equation}
v_{n}\geq \left( 1+\frac{k}{n}\right)n^{a} \label{vlower}
\end{equation}
for $n=N-1$ and $n=N-2$, where $a=-1+1/\eps$ and $k=1+\lam/\eps$,
then the same inequality holds for all $n\geq N$.
\end{lemma}

\Proof Suppose that $n\geq \lam+3$ and that
\[
v_{n-i}\geq \left( 1+\frac{k}{n-i}\right) (n-i)^a
\]
for $i=1,\, 2$. Then
\begin{eqnarray*}
n^{-a}v_n & = & n^{-a}\left( \frac{n-2}{n}
v_{n-2}+2\frac{n-1-\lam}{\eps
n(n-1)}v_{n-1}\right)\\
&\geq & \left(1+\frac{k}{n-2}\right) \left(1-\frac{2}{n}\right)^{a+1}%
+2\frac{n-1-\lam}{\eps n(n-1)}\left(1+\frac{k}{n-1}\right)\left(1-\frac{1}{n}\right)^a\\
&= & 1+\frac{k}{n}+\frac{2}{n^2}+O(n^{-3})\\
&\geq& 1 +\frac{k}{n}
\end{eqnarray*}
provided $n\geq N=N_{\lam,\eps}^{(2)}$.  It follows inductively that
(\ref{vlower}) holds for all $n\geq N$.

\begin{theorem}\label{Apluspositive} If $\, 0<\eps<2$ and
$\lam$ is a real eigenvalue of $A_+$ then $\lam > 1$.
\end{theorem}

\Proof Suppose that $A_+v=\lam v$ where $\lam\leq 1$ and $v_1=1$.
The initial condition $\eps v_2=(1-\lam)v_1$ implies that $v_2\geq
0$, and it then follows from the signs of the coefficients in
(\ref{recurrence}) that $v_n>0$ for all $n\geq 3$.
Lemma~\ref{vlower} implies that there exists a constant $c>0$ such
that
\[
v_{n}\geq c\left( 1+\frac{k}{n}\right)n^{a}
\]
for all $n\geq N_{\lam,\eps}^{(2)}$. The lower bound $a>-1/2$
implies that $v\notin l^2(\Z_+)$, and hence that $\lam$ is not an
eigenvalue of $A_+$.

{\bf Hypothesis} From this point onwards we assume that $0<\eps<2$
and $\lam\geq 0$.

\begin{theorem}\label{GrowAsymp} For every $\del >0$ there exists
$N=N_{\lam,\eps,\del}$ and a solution $v$ of (\ref{recurrence}) such
that
\[
n^a\leq v_n\leq (1+\del)n^a
\]
for all $n\geq N$, where $a=-1+1/\eps$.
\end{theorem}

\Proof We put $N=N_{\lam,\eps,\del}=\max%
\{ N^{(1)}_{\lam,\eps},N^{(2)}_{\lam,\eps},2+k/\del\}$ where
$k=1+\lam/\eps$ and let $v$ be the solution of (\ref{recurrence})
such that $v_{N-i}=(1+\del)(N-i)^a$ for $i=1,\, 2$.
Lemma~\ref{GrowUpper} implies that $0<v_n\leq (1+\del)n^a$ for all
$n\geq N$. Since
\[
v_{n}\geq \left( 1+\frac{k}{n}\right)n^a
\]
for $n=N-1$ and $n=N-2$, we deduce by Lemma~\ref{GrowLower} that
$v_n\geq (1+k/n)n^a\geq n^a$ for all $n\geq N$. This completes the
proof.

We will show that, up to a multiplicative constant, there is exactly
one `subordinate' solution $v$ of (\ref{recurrence}) such that
$\lim_{n\to +\infty} v_n=0$. We identify this solution by solving
the recurrence relation backwards from $n=M$ and then letting
$M\to+\infty$.

\begin{lemma}\label{SuborUpper} There exists $N=N_{\lam,\eps}^{(3)}$
such that if $M>N$ and $v_n=(-1)^nw_n$ is a solution of
(\ref{recurrence}) satisfying $0<w_{M+i}\leq(M+i)^{-c}$ for $i=1,\,
2$ where $c=1+1/\eps$, then $0< w_n\leq n^{-c}$ for all $n$
satisfying $N\leq n\leq M$.
\end{lemma}

\Proof The sequence $w_n$ satisfies the recurrence relation
\begin{equation}
w_n=\frac{n+2}{n}w_{n+2}%
+\frac{2(n+1-\lam)}{\eps n(n+1)} w_{n+1}\, .\label{wrecurrence}
\end{equation}
This has positive coefficients for $n\geq \lam$ so the solution is
positive if $\lam<n\leq M$. Suppose inductively that $0<w_{n+2}\leq
(n+2)^{-c}$ and $0<w_{n+1}\leq (n+1)^{-c}$ for such an $n$. Then
\begin{eqnarray*}
n^cw_n & \leq &
\left(1+\frac{2}{n}\right)^{1-c}+\frac{2(n+1-\lam)}{\eps
n(n+1)}\left(1+\frac{1}{n}\right)^{-c}\\
&=& 1-\frac{2\lam}{\eps n^2}+O(n^{-3})\\
&\leq & 1
\end{eqnarray*}

for all large enough $n$. By induction there exists
$N=N_{\lam,\eps}$ such that $0<w_n\leq n^{-c}$ provided $N\leq n\leq
M$.

\begin{lemma}\label{SuborLower} There exists $N=N_{\lam,\eps}^{(4)}$
such that if $M>N$ and $v_n=(-1)^nw_n$ is a solution of
(\ref{recurrence}) such that
\begin{equation}
w_n\geq \left( 1-\frac{h}{n}\right) n^{-c} \label{wlower}
\end{equation}
for $n=M+1$ and $n=M+2$, where $c=1+1/\eps$ and $h=1+\lam/\eps$,
then (\ref{wlower}) holds for all $n$ satisfying $N\leq n\leq M$.
\end{lemma}

\Proof Suppose that $\max\{h,\lam\}\leq n\leq M$ and (\ref{wlower})
holds when $n$ is replaced by $n+1$ or $n+2$. Then
\begin{eqnarray*}
n^cw_n&\geq& \left(1-\frac{h}{n+2}\right) \left(
1+\frac{2}{n}\right)^{1-c}+2\frac{n+1-\lam}{\eps n(n+1)}
\left(1-\frac{h}{n+1}\right)\left(1+\frac{1}{n}\right)^{-c}\\
&=& 1-\frac{h}{n}+\frac{2}{n^2}+O(n^{-3})\\
&\geq& 1-\frac{h}{n}
\end{eqnarray*}
provided $n$ is large enough. An induction now implies that there
exists $N=N_{\lam,\eps}^{(4)}$ such that (\ref{wlower}) holds for
all $n$ such that $N\leq n\leq M$.

\begin{theorem}\label{subordinate}
There exists $N=N_{\lam,\eps}^{(5)}$ and a unique solution
$v_n=(-1)^n w_n$ of (\ref{recurrence}) such that
\[
\left( 1-\frac{h}{n}\right)n^{-c} \leq w_n\leq n^{-c}
\]
for all $n\geq N$, where $c=1+1/\eps$ and $h=1+\lam/\eps$. Hence
\begin{equation}
\lim_{n\to +\infty} w_nn^c=1. \label{normalize}
\end{equation}
\end{theorem}

\Proof Let $M
> N=N_{\lam,\eps}^{(5)}
=\max\{N_{\lam,\eps}^{(3)},N_{\lam,\eps}^{(4)}\}$ and let $w^{(M)}$
denote the solution of (\ref{wrecurrence}) such that
$w^{(M)}_n=n^{-c}$ for $n=M+1$ and $n=M+2$. Lemmas~\ref{SuborUpper}
and \ref{SuborLower} imply that
\[
\left( 1-\frac{h}{n}\right)n^{-c} \leq w_n^{(M)}\leq n^{-c}
\]
for all $n$ such that $N\leq n\leq M$. By choosing a sequence
$M_r\to +\infty$ such that $w_{N}^{(M_r)}$ and $w_{N+1}^{(M_r)}$
converge as $r\to +\infty$ we see using (\ref{wrecurrence}) that
$w_{n}^{(M_r)}$ converge for all $n\geq 1$. Denoting the limit by
$w^{(\infty)}$ we deduce that
\[
\left( 1-\frac{h}{n}\right)n^{-c} \leq w^{(\infty)}_n\leq n^{-c}
\]
for all $n\geq N$. Putting $v^{(\infty)}_n=(-1)^nw^{(\infty)}_n$,
the uniqueness of the solution $v^{(\infty)}$ subject to the
normalization condition (\ref{normalize}) follows from the fact that
the solution space of (\ref{recurrence}) is two-dimensional and it
contains a divergent sequence by Theorem~\ref{GrowAsymp}.

Numerical examples suggest that the following lemma is not the best
possible and that $w$ takes its maximum value very close to
$n=\lam$. Figure~1 plots the eigenfunction $v$ of the operator $A_+$
for the eigenvalue $\lam\sim 14.94784$ with $\eps=0.1$.

\begin{center}
\scalebox{0.5}{\includegraphics{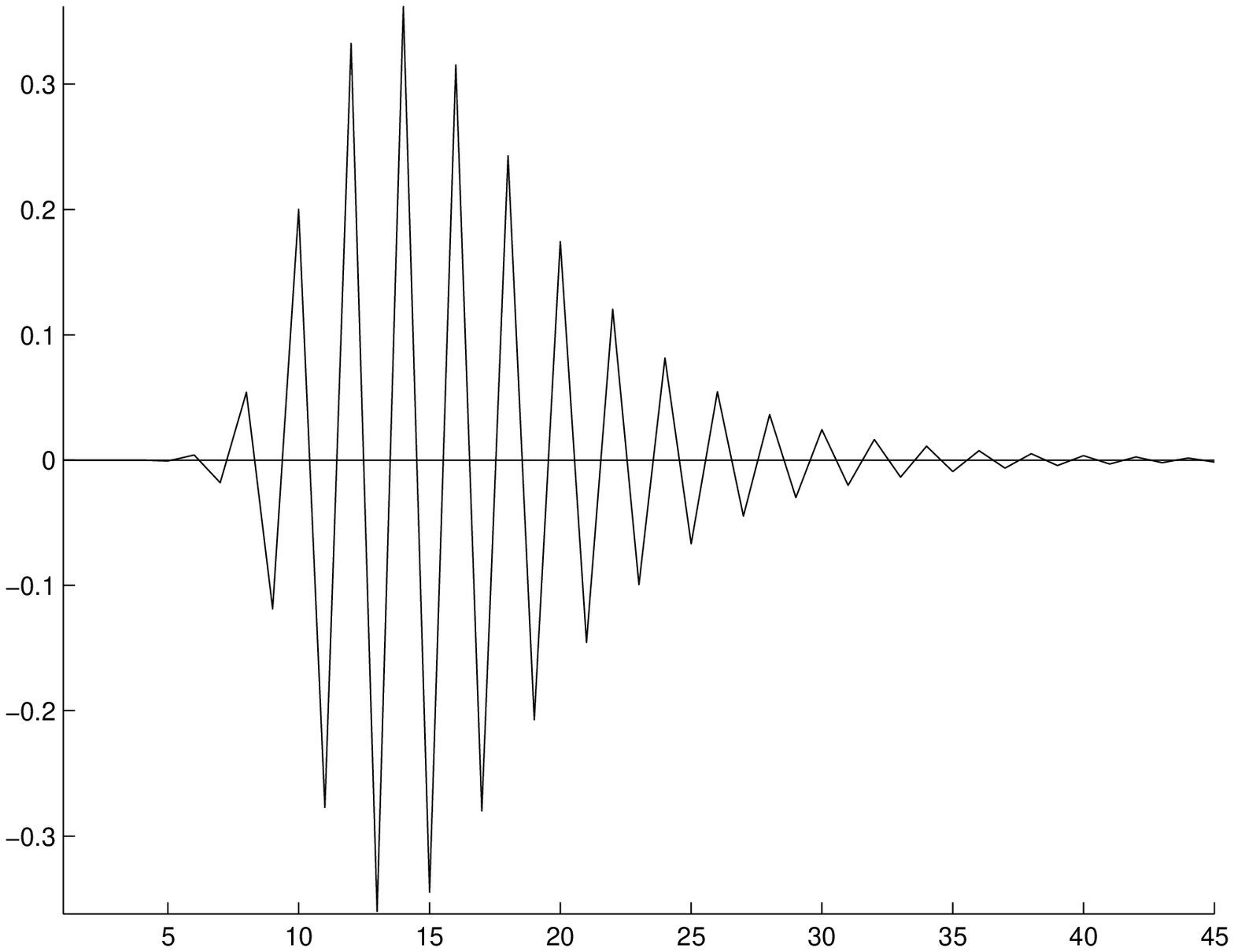}}\\
 Figure~1. Eigenvector $v$ for $\eps=0.1$ and $
\lam\sim 14.94784$   % uses sazonov24.m
\end{center}

\begin{lemma}\label{monotone}
If $\lam\geq 0$ then the unique subordinate solution $w$ of
(\ref{wrecurrence}) satisfies $ 0<w_{n+1}<w_n $ for all $n\geq
2\lam$.
\end{lemma}

\Proof Let $w^{(M)}$ denote the solution of (\ref{wrecurrence})
constructed in the proof of Theorem~\ref{subordinate}. Then
\begin{eqnarray*}
w^{(M)}_M&=& \frac{M+2}{M}(M+2)^{-c}+\frac{2(M+1-\lam)}{\eps
M(M+1)}(M+1)^{-c}\\
&=& (M+1)^{-c}\left( (1+2/M)^{1-c}(1+1/M)^c+\frac{2}{\eps
M} (1-\lam/(M+1)) \right)\\
&=& (M+1)^{-c}\left( 1+c/M+O(M^{-2})\right)\\
&\geq& (M+1)^{-c}
\end{eqnarray*}
provided $M$ is large enough. Therefore \[w_M^{(M)}\geq
w_{M+1}^{(M)}\geq w_{M+2}^{(M)}.\]

We prove inductively that $w_n^{(M)}\geq w_{n+1}^{(M)}$ for all $n$
such that $2\lam\leq n\leq M$. If this holds with $n$ replaced by
$n+1$ or by $n+2$ then
\begin{eqnarray*}
w_n^{(M)}-w_{n+1}^{(M)} &=& \frac{n+2}{n}w_{n+2}^{(M)}
+\frac{2(n+1-\lam)}{\eps n(n+1)}w_{n+1}^{(M)}\\
&& \, - \frac{n+3}{n+1}w_{n+3}^{(M)}
-\frac{2(n+2-\lam)}{\eps (n+1)(n+2)}w_{n+2}^{(M)}\\
&=&\frac{n+2}{n}w_{n+2}^{(M)}-\frac{n+3}{n+1}w_{n+3}^{(M)}\\
&&\, +\frac{2(n+1-\lam)}{\eps n(n+1)}w_{n+1}^{(M)}-%
\frac{2(n+2-\lam)}{\eps (n+1)(n+2)}w_{n+2}^{(M)}\\
&\geq& \left( \frac{n+2}{n}-\frac{n+3}{n+1}\right)
w_{n+3}^{(M)}\\
&&\, +\left( \frac{2(n+1-\lam)}{\eps n(n+1)}-%
\frac{2(n+2-\lam)}{\eps (n+1)(n+2)}\right) w_{n+2}^{(M)}\\
&\geq 0
\end{eqnarray*}
provided $n\geq 2\lam$. This completes the induction.

Finally we take the same sequence $M(r)$ as in the proof of
Theorem~\ref{subordinate} to obtain $0<w_{n+1}^{(\infty)}\leq
w_n^{(\infty)}$ for all $n\geq 2\lam$.

%%%%%%%%%%%%%%%%%%%%%%%%%%%%%%%%%%%%%%%%%%%%%%%%%%%%%%%%%%%%%%%%%%
\section{Compactness of the
Resolvent\label{SectResolvent}}

In this section we prove that $0\notin\Spec(A_+)$ and that
$A_+^{-1}$ is a Hilbert-Schmidt operator, and hence compact. This
implies that the spectrum of $A_+$ is discrete and coincides with
its set of eigenvalues. We cannot, however, prove that the spectrum
is real. We define the Hilbert-Schmidt operator $R$ on $l^2(\Z_+)$
by
\begin{equation}
(Rf)_m=\sum_{n=1}^\infty \rho_{m,n}f_n\label{Rdef}
\end{equation}
where $\rho\in l^2(\Z_+\times \Z_+)$ is given explicitly. We then
show directly that $R$ is the inverse of $A_+$.

Let $\phi$ be the solution of
\[
\phi_n=\frac{n-2}{n}\phi_{n-2}+\frac{2}{\eps n}\phi_{n-1}
\]
that satisfies the initial conditions $\phi_1=1$ and
$\phi_2=\eps^{-1}$. One sees immediately that $\phi_n>0$ for all
$n\geq 1$. Theorem~\ref{GrowAsymp} implies that there exists a
constant $c_1>0$ such that
\[
c_1^{-1}n^a\leq \phi_n\leq c_1 n^a
\]
for all $n\geq 1$.

Let $\psi_n=(-1)^n w_n$ be the unique subordinate solution of
\[
\psi_n=\frac{n-2}{n}\psi_{n-2}+\frac{2}{\eps n}\psi_{n-1}
\]
such that $w$ satisfies the asymptotic condition $\lim_{n\to
+\infty}n^cw_n=1$. Since
\[
w_n= \frac{n+2}{n}w_{n+2}+\frac{2}{\eps n}w_{n+1}
\]
we see that $w_n>0$ for all $n\geq 1$, and indeed that there exists
a constant $c_2>0$ such that
\[
c_2^{-1}n^{-c}\leq w_n\leq c_2 n^{-c}
\]
for all $n\geq 1$.

We finally put
\[
\sig_n=\frac{\eps}{2}n(n-1)\phi_{n-1}w_n+\frac{\eps}{2}n(n+1)\phi_nw_{n+1}+n\phi_nw_n
\]
and observe that $\sig_n>0$ for all $n\geq 1$. The upper and lower
bounds on $\phi$ and $w$ imply that there exists a constant $c_3>0$
such that
\[
c_3^{-1}n^{a-c+2}\leq \sig_n\leq c_3 n^{a-c+2}
\]
for all $n\geq 1$.

\begin{theorem}\label{resolventbounds}
If $0<\eps<2$ and
\[
\rho_{m,n}=\left\{ \begin{array}{ll}%
(-1)^n\phi_m\psi_n/\sig_n & \mbox{ if $m\leq n$,}\\
(-1)^n\psi_m\phi_n/\sig_n & \mbox{ if $m> n$.}
\end{array}  \right.
\]
then $\rho\in l^2(\Z_+\times \Z_+)$. The Hilbert-Schmidt operator
$R$ defined by (\ref{Rdef}) satisfies $A_+ Rf=f$ for all $f\in
l^2(\Z_+)$. Indeed $0\notin \Spec(A_+ )$ and $R=A_+^{-1}$.
\end{theorem}

\Proof The above bounds on $\phi,\, \psi,\, \sig$ imply that
\[
|\rho_{m,n}|\leq\left\{ \begin{array}{ll}%
c_4 m^an^{-a-2} & \mbox{ if $m\leq n$,}\\
c_4 m^{-c}n^{c-2}& \mbox{ if $m> n$.}
\end{array} \right.
\]
It follows that
\[
\sum_{m=1}^\infty |\rho_{m,n}|^2\leq c_5 n^{-3}
\]
and then that
\[
\sum_{m,n=1}^\infty |\rho_{m,n}|^2 <\infty.
\]
We conclude that $R$ is a compact operator. If
$\{e_n\}_{n=1}^\infty$ is the standard basis in $l^2(\Z_+)$ then a
direct calculation shows that $A_+ Re_n=e_n$ for all $n$. By using
the fact that $A_+ $ is closed one deduces that $\Ran(R)\subseteq
\Dom(A_+)$ and that $A_+ Rf=f$ for all $f\in l^2(\Z_+)$. We conclude
from this that $\Ran(A_+ )=l^2(\Z_+)$. The bound $0<\eps<2$ implies
that $\Ker(A_+ )=\{ 0\}$ by Theorem~\ref{Apluspositive}, so we
finally see that $0\notin \Spec(A_+ )$ and that $R=A_+^{-1}$.

%%%%%%%%%%%%%%%%%%%%%%%%%%%%%%%%%%%%%%%%%%%%%%%%%%%%%%%%%%%%%%
\section{$\lam$-Dependence}
\newcommand{\wt}{\tilde{w}}%
\newcommand{\lamt}{\tilde{\lambda}}%

In this section we prove that the unique normalized subordinate
solution $v_{\lam,n}=(-1)^n w_{\lam,n}$ of (\ref{recurrence})
provided by Theorem~\ref{subordinate} depends continuously on
$\lam$.

We first observe that for any $\Lam\geq 1$ the various constants
$N_{\lam,\eps}^{(i)}$ are uniformly bounded with respect to $\lam$
provided $0\leq \lam\leq \Lam$. We (incorrectly) use the notation
$N_{\Lam,\eps}^{(i)}$ to refer to the relevant upper bounds.

\begin{lemma}\label{cont1}
If $0\leq\lam\leq\mu\leq \Lam$ then
\[
0<w_{\Lam,n}\leq w_{\mu,n}\leq w_{\lam,n}\leq w_{0,n}<\infty
\]
for all $n\geq \Lam$.
\end{lemma}

\Proof The positivity of $w_{\lam,n}$ for $n\geq \Lam$ follows from
the positivity of the coefficients of (\ref{wrecurrence}) for $n\geq
\Lam$ and the positivity of $w_{\lam,n}$ for all $n\geq
N=N_{\Lam,\eps}^{(5)}$. We only need only prove the central
inequality above since the other two are special cases of it.

Theorem \ref{subordinate} implies that if $\del >0$ then
\begin{equation}
w_{\mu,n}\leq (1+\del)w_{\lam,n} \label{almostmono}
\end{equation}
for all $n\geq N=N_{\Lam,\eps,\del}^{(6)}$. This inequality persists
for all $n\in [\Lam,N]$ by the monotonicity of the coefficients of
(\ref{wrecurrence}). Since (\ref{almostmono}) holds for all $\del
>0$ and all $n\geq \Lam$, the required inequality follows by letting
$\del\to 0$.

\begin{lemma}\label{cont2}
If $\, 0\leq\lam\leq\mu\leq \Lam$ and $|\mu- \lam|\leq \del$ then
\begin{equation}
0<w_{\lam,n}\leq p_{\Lam,\eps,n,\del}w_{\mu,n}\label{almostupper}
\end{equation}
for all $n\geq 2\Lam$, where
\[
p_{\Lam,\eps,n,\del}=(1+\del)\exp\left\{2\del\eps^{-1}
\sum_{r=n}^\infty r^{-2}\right\}.
\]
\end{lemma}

\Proof Since $1+\del\leq p_{\Lam,\eps,n,\del}$, Theorem
\ref{subordinate} implies that (\ref{almostupper}) holds for all
$n\geq N=N_{\Lam,\eps,\del}^{(6)}$. We prove inductively that the
same inequality persists for $n\in [2\Lam,N]$. If
(\ref{almostupper}) holds with $n$ replaced by $n+1$ and by $n+2$,
then, using Lemma~\ref{monotone}, we obtain
\begin{eqnarray*}
w_{\lam,n}&=&\frac{n+2}{n}w_{\lam,n+2}+\frac{2(n+1-\lam)}{\eps
n(n+1)}w_{\lam,n+1}\\
&\leq &
\frac{n+2}{n}p_{\Lam,\eps,n+2,\del}w_{\mu,n+2}+\frac{2(n+1-\lam)}{\eps
n(n+1)}p_{\Lam,\eps,n+1,\del}w_{\mu,n+1}\\
&\leq &
p_{\Lam,\eps,n+1,\del}\left(\frac{n+2}{n}w_{\mu,n+2}+\frac{2(n+1-\lam)}{\eps
n(n+1)}w_{\mu,n+1}\right)\\
&\leq &
p_{\Lam,\eps,n+1,\del}\left(\frac{n+2}{n}w_{\mu,n+2}+\frac{2(n+1-\mu)}{\eps
n(n+1)}w_{\mu,n+1}+ \frac{2\del}{\eps n^2} w_{\mu,n+1}\right)\\
&\leq &
p_{\Lam,\eps,n+1,\del}\left(w_{\mu,n}+\frac{2\del}{\eps n^2} w_{\mu,n+1}\right)\\
&\leq & p_{\Lam,\eps,n+1,\del}\left(1+\frac{2\del}{\eps
n^2}\right)w_{\mu,n}\\
&\leq & p_{\Lam,\eps,n,\del}w_{\mu,n}.
\end{eqnarray*}
This completes the induction.

\begin{theorem}\label{cont3} The subordinate solution
$v_\lam$ depends continuously on $\lam$ for $0\leq \lam <\infty$.
Hence the function
\[
f(\lam):=\eps v_{\lam,2}-(1-\lam) v_{\lam,1}
\]
is continuous on $[0,\infty)$.
\end{theorem}

\Proof It is sufficient to prove that $f$ is continuous on
$[0,\Lam]$ for every positive integer $\Lam$. It follows directly
from the estimates in Lemmas~\ref{cont1} and \ref{cont2} that the
map $\lam\in [0,\Lam]\to (w_{\lam,2\Lam},w_{\lam,2\Lam+1})$ is
continuous. Composing this with the linear (and therefore
continuous) map $(w_{\lam,2\Lam},w_{\lam,2\Lam+1})\to \eps
v_{\lam,2}-(1-\lam) v_{\lam,1}$ yields the second statement of the
theorem.

%%%%%%%%%%%%%%%%%%%%%%%%%%%%%%%%%%%%%%%%%%%%%%%%%%%%%
\section{Numerical Calculations}

Let $v_\lam$ denote the solution of (\ref{recurrence}) such that
\[
\lim_{n\to+\infty} v_{\lam,n}(-1)^n n^c=1.
\]
Then $\lam>0$ is an eigenvalue if and only if
\[
f(\lam):= \eps v_{\lam,2}-(1-\lam) v_{\lam,1}
\]
vanishes. Since this function is continuous, one can compute the
roots of $f(\lam)=0$ by evaluating $f(\lam)$ numerically for a range
of values of $\lam$. We determined the subordinate solution by
solving (\ref{wrecurrence}), starting from $M=4000$ (and also
$M=8000$ to check consistency) with $w_{M+i}=(M+i)^{-c}$ for $i=1,\,
2$. Figure~2 plots $f(\lam)$ for $\eps=0.1$ and $0\leq \lam\leq 4$.
The eigenvalues listed in Table 1 were obtained by solving
$f(\lam)=0$ numerically, and are quite close to those obtained in
\cite{BOS}.

\begin{center}
\begin{tabular}{ccc}
$n$&$\lam_n$& $\norm P_n\norm$\\
1&1.00968&1.0189\\
2&2.07334&1.1848\\
3&3.22978&1.8868\\
4&4.50134&4.3409\\
5&5.89993&13.341\\
6&7.43194&50.638\\
7&9.10097&226.20\\
8&10.9092&1152.9\\
9&12.8578&6561.3\\
10&14.9478&41018\\
15&27.5331&$-$\\
20&43.74&$-$
\end{tabular}\\
{ Table 1. Eigenvalues of $A_+$ for $\eps=0.1$}  % uses sazonov17.m
\end{center}

The computation is very stable and one can confidently evaluate the
first ten eigenvalues to much higher accuracy. The list of
eigenvalues found is compatible with the asymptotic formula $
\lam_n\sim \alp n^{\gam} $ where $\alp\sim0.53$ and $\gam\sim 1.44$.

However, for $\eps=1$, the Fourier coefficients decrease much more
slowly, and the eigenvalue calculation is correspondingly more
onerous. We computed the first five eigenvalues for $\eps=1$,
determining the subordinate solution as before with $M$ between
$1000$ and $32000$. The apparent numbers of eigenvalues increased
from $7$ to $11$ as $M$ increased in this range. For $M=4000$ it
appeared that the computation of the first five eigenvalues
presented in Table 2 was reliable.

\begin{center}
\begin{tabular}{cc}
$n$&$\lam_n$\\
1&1.4485\\
2&4.3159\\
3&8.6219\\
4&14.3638\\
5&21.5414
\end{tabular}\\
{ Table 2. Eigenvalues of $A_+$ for $\eps=1$}  % uses sazonov21.m
\end{center}

\begin{center}
\scalebox{0.6}{\includegraphics{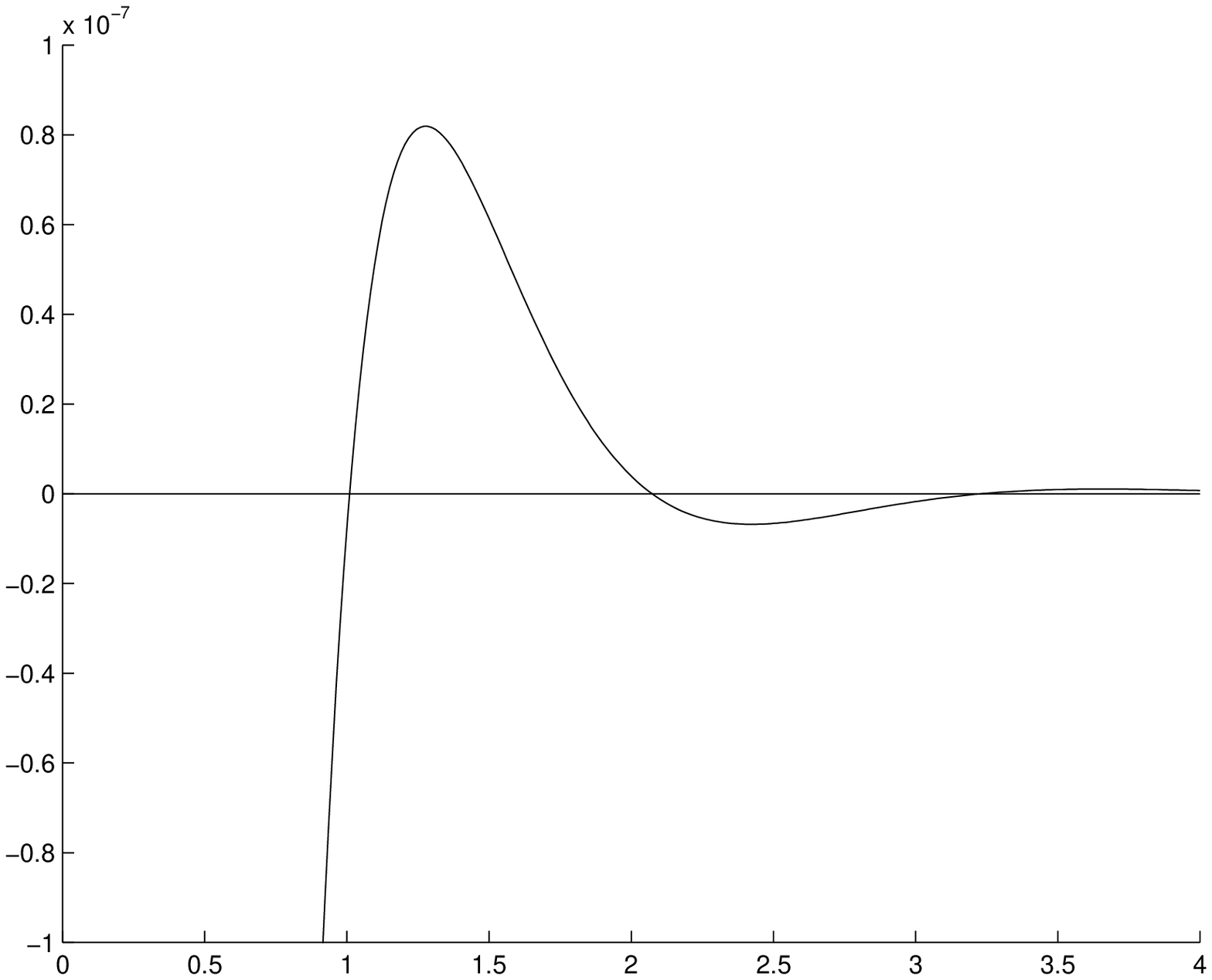}}\\
 Figure~2. $f(\lam)$ for $\eps=0.1$ and $0\leq
\lam\leq 4$   % uses sazonov12.m
\end{center}

We conclude with some comments about the conjecture in \cite{BOS}
that the eigenvectors form a basis. It seems quite plausible that
they form a complete set in the sense that their linear span is
dense. However, if they form a basis then the spectral projections
\[
P_nf=\frac{\langle f,\phi^\ast_n\rangle} {\langle
\phi_n,\phi^\ast_n\rangle}\phi_n
\]
of $H$ must be uniformly bounded in norm, where $\phi_n$ are the
eigenfunctions of $H$ and $\phi_n^\ast$ the corresponding
eigenfunctions of $H^\ast$; see \cite[Lemma~3.3.3]{LOTS}. However,
it appears from \cite[Figure~4]{BOS} that the eigenfunctions
$\phi_n$ concentrate more and more strongly around $\theta=\pi$ as
$n$ increases; the eigenfunctions $\phi_n^\ast$ should concentrate
around $\theta=0$ as $n\to\infty$ for similar reasons. If this is
indeed the case then the norms of the spectral projections
\[
\norm P_n\norm =\frac{\norm \phi_n\norm \, \norm \phi_n^\ast\norm}%
{|\langle \phi_n,\phi^\ast_n\rangle|}
\]
must diverge as $n\to\infty$ and the eigenfunctions do not form a
basis.  The norms of the first $10$ spectral projections are
presented in Table~1 and confirm the conjecture that they diverge as
$n$ increases. See \cite{DK} for another highly non-self-adjoint
operator arising in physics for which an apparently well-behaved
sequence of eigenfunctions do not form a basis.

\vspace{4mm} {\bf Acknowledgements} I should like to thank I A
Sazonov for drawing my attention to the results in \cite{BOS}.

Department of Mathematics\\
King's College\\
Strand\\
London WC2R 2LS

E.Brian.Davies@kcl.ac.uk
\end{document}